\documentclass[11pt]{article}
\usepackage{amsfonts}
\usepackage{amssymb}
\usepackage{amscd}
\usepackage{amsmath}
\usepackage{epsfig}
\usepackage{graphicx, subfigure}
\usepackage{latexsym}
\usepackage{pst-3dplot}
\usepackage[symbol*]{footmisc}

\usepackage{color}
\usepackage{xcolor}
\usepackage{epsfig}
\usepackage{graphicx, subfigure}
\usepackage[makeroom]{cancel}
\usepackage{ulem}

\newlength{\dinwidth}
\newlength{\dinmargin}
\usepackage[all,cmtip]{xy}

\newtheorem{definition}{Definition}
\newtheorem{theorem}{Theorem}

\newtheorem{corollary}{Corollary}

\newtheorem{lemma}{Lemma}
\newtheorem{example}{Example}

\def\Z{\mathcal Z}

\def \P{\mathbb P}

\def\M{\mathcal M}
\def\D{\mathcal V}

\begin{document}

\renewcommand*{\thefootnote}{\fnsymbol{footnote}}

\title{Mutual incidence matrix of two balanced incomplete block designs}
\author{A. Shramchenko, V. Shramchenko$^*$}

\date{}

\maketitle

\footnotetext[1]{Department of mathematics, University of
Sherbrooke, 2500, boul. de l'Universit\'e,  J1K 2R1 Sherbrooke, Quebec, Canada. E-mail: {\tt Vasilisa.Shramchenko@Usherbrooke.ca}}

 \begin{abstract}
We propose to consider a mutual incidence matrix $M$ of two balanced incomplete block designs built on the same finite set. In the simplest case, this matrix reduces to the standard incidence matrix of one block design. We find all eigenvalues of the  matrices $MM^T$ and $M^TM$ and their eigenspaces. 
 
  \end{abstract}
\maketitle

\section{Introduction}

Incidence matrices are an important tool in graph theory and other mathematical contexts. In the theory of  balanced incomplete block designs (BIBDs), see \cite{encyclopedia, handbook, Rib} for a detailed introduction, the incidence matrix is a matrix whose rows correspond to points and columns correspond to blocks. The entry in the row $i$ and column $j$ is 1 if the point $i$ belongs to the block $j$ and 0 otherwise. 
Such incidence matrices can be used, for example,  to prove Fisher's inequality, a fundamental theorem giving a necessary condition for the existence of a balanced incomplete block design, stating that the number of blocks is greater or equal to the number of points, see for example \cite{R}. A number of results on existence of block designs is obtained using incidence matrices, see Section V.6.4 of \cite{handbook}. A slightly modified notion of an incidence matrix is used in the context of subspace designs, see \cite{subspace}.

The aim of this paper is to introduce and study  a {\it mutual incidence matrix}  $M(D_1, D_2)$ for two block designs $D_1,\, D_2$ built on the same set $V$ of finite cardinality $v$. The rows of the matrix  are indexed by blocks of $D_1$ and the columns  are indexed by blocks of $D_2$. The entry on the intersection of a row and a column is the number of elements shared by the corresponding blocks of the two block designs. Interchanging the two block designs amounts to the transposition of the mutual incidence matrix.

 Naturally, for the simplest block design $\D$ having one-element blocks where each element of the set $V$ is a block, the mutual incidence matrix $M(\D, D)$ is the incidence matrix of the block design $D$. 

The mutual incidence matrix $M=M(D_1,D_2)$ not being a square matrix in general, we consider  symmetric matrices $MM^T$ and $M^TM$. The paper is devoted to the study of  these matrices: we find all  their eigenvalues and eigenspaces in terms of the parameters of the two block designs, see Theorem \ref{thm}.  In particular, using properties of block designs, we construct two sets of vectors $\Z_{D_1}(x,y)$ and $\Z_{D_2}(x,y)$ for $x,y\in V$, see \eqref{Z1} and \eqref{Z1_comp}.
We show that the vectors $\Z_{D_1}(x,y)$ are eigenvectors of the matrix $MM^T$ and vectors $\Z_{D_2}(x,y)$ are eigenvectors of $M^TM$ and find the corresponding eigenvalues in terms of the parameters of the two block designs $D_1$ and $D_2$. Note that these results apply to matrices $MM^T$ and $M^TM$ with $M$ being a mutual incidence matrix of two block designs constructed on the same set $V$ without assuming any particular ordering of the blocks of the two block designs.

For the incidence matrix $M=M(\D, D)$ of a block design $D$ built on the set $V$, the symmetric matrix $MM^T$ has a simple structure and is of rank $v.$  Due to this simple structure of $MM^T$, its eigenvalues and the corresponding eigenvectors are easily obtained \cite[\S VII.7]{handbook}. On the other hand, the structure of the other symmetric matrix $M^TM$ obtained from $M$ is not as simple. While the two matrices, $MM^T$ and $M^TM$ have the same nonzero eigenvalues, finding the corresponding eigenvectors of $M^TM$ is not straightforward and, to our knowledge, has not been addressed. We obtain them in the present work, see Corollary \ref{corollary}.

Given the importance of the incidence matrix of a block design in studying designs, one may expect that mutual incidence matrices may also be used to understand properties of block designs. The outlined results of the paper are quite general, that is they apply to any two block designs built on the same set. The fact that one can obtain such a general result is our main motivation for introducing the mutual incidence matrix as it points to the fact that such matrices might be of particular interest and merit further study.

The incidence matrix is closely related to the incidence graph of a block design, see for example, \cite{crn,nato, ma, domination, igraph} and references therein for definitions and related results. This, in particular, allows to apply results and terminology from graph theory to the theory of block designs.
While entries of the incidence matrix for any block design are equal to $0$ and $1$, the entries of the mutual incidence matrix can be any non-negative integers. Thus it corresponds to a multigraph, the {\it mutual incidence graph} of two block designs, that is a graph having multiple edges. 

\section{Block designs and their mutual incidence matrix and graph}
\label{sect_def}

A balanced incomplete block design is defined as follows, see \cite{encyclopedia, handbook, Rib}.

\begin{definition}
\label{def}
Let $V$ be a finite set of cardinality $v=|V|$ and let $b,k,r,\lambda$ be four positive integers with $k<v$. A (balanced incomplete) block design with parameters $(v,b,r,k,\lambda)$ built on the set $V$ is a list of $b$ blocks, each of which is a $k$-element subset of $V$, such that every element of $V$ is contained in exactly $r$ blocks and every pair of elements of $V$ is contained in exactly $\lambda$ blocks.
\end{definition}

Note that the parameters $v,b,r,k,\lambda$ of a given block design are not independent;  the following relations follow from Definition \ref{def}: $bk=vr$ and $r(k-1)=\lambda(v-1)$. Thus it suffices to specify parameters $v, k, \lambda$ for a block design. The structure from Definition \ref{def} is sometimes referred to as $2-(v,k,\lambda)$ design, or simply a 2-design, or a pairwise balanced design, and is a particular case of a $t-(v,k,\lambda)$ design where $\lambda$ is the number of times any given set of $t$ elements (points) of $V$ appears in the blocks of the design.  We consider the set $V$ itself as a simplest case of a block design with $b=v,\;$ $r=k=1$ and $\lambda=0\,.$
\\

We denote $D_i, i \in \mathbb N,$ a block design from Definition \ref{def} built on the set $V$ with parameters  $(v,b_i,r_i, k_i,\lambda_i)$. The blocks of $D_i$ are assumed ordered in some arbitrary way from $1$ to $b_i$ and denoted by $B^s_i$ with  $s=1, \dots,  b_i.$

For two block designs $D_1$ and $D_2$ built on the same finite set $V$ define the {\it mutual incidence matrix} $M=M(D_1, D_2)$ as a matrix of the size $b_1\times b_2$ whose rows are indexed by the set of blocks of $D_1$ and columns are indexed by the set of blocks of $D_2$ and such that the entry $M_{ij}$ in the row corresponding to $B_1^i$ and the column corresponding to $B_2^j$ is
\begin{equation}
\label{Mij}
M_{ij}(D_1, D_2)=M_{ij}:= |B_1^i\cap B_2^j|,
\end{equation}
the number of common elements in the two blocks.

Our objective is to find eigenvalues of the square matrix $MM^T$ obtained from the above mutual incidence matrix of two BIBDs. Let us first formulate a simple but useful lemma. 
\begin{lemma}
\label{lemma_main}
Let $Q=\{Q_1, \dots, Q_N\}$ be some set of subsets of $V$. For $x\in V$ define a function $t_Q(x)$  by setting $t_Q(x)$ to be the number of elements of $Q$ containing $x$:
\begin{equation}
\label{t-def}
t_Q(x) = \sum_{j=1}^N |Q_j\cap x|\,.
\end{equation}
 Then  for any subset $X\subseteq V$
\begin{equation}
\label{t}
\sum_{j=1}^N |Q_j\cap X|=\sum_{x\in X}t_Q(x)\,.
\end{equation}

\end{lemma}

Note that the classical incidence matrix $M$ of any block design with parameters $(v,b,r,k,\lambda)$ satisfies $MM^T=(r-\lambda)I +\lambda J$ where $I$ is the $v\times v$ identity matrix and $J$ is the matrix whose entries are all equal to $1$, see \cite{handbook}. In particular, all diagonal entries of this matrix are equal to $r$ and all off-diagonal entries are equal to $\lambda.$ 

The next theorem shows that this is partly true for the mutual incidence matrix $M=M(D_1, D_2)$ of two block designs as well: all diagonal entries of the matrix $MM^T$ are equal and can be expressed in terms of the parameters of $D_1$ and $D_2$. The off-diagonal entries however may take more than one value as illustrated in Example \ref{example3} in Section \ref{sect_examples}.
Let us start with a lemma. 
\begin{lemma} \label{lemma_formulas}
 Let $D$ be a block design built on a finite set $V$ with parameters $(v,b,r,k,\lambda)$ and let $S\subseteq V$ be any subset of $V$ of cardinality $|S|$. For $0\leq j\leq k$ let $\delta_j(S, D)$ be the number of blocks of $D$ sharing exactly $j$ elements with $S$. Then
\begin{equation*}
\sum_{j=0}^k {j}\delta_j(S, D)=r|S|\qquad\text{ and }\qquad
\sum_{j=0}^k {j^2}\delta_j(S, D)=(\lambda{|S|}-\lambda+r)|S|.
\end{equation*} 
\end{lemma}
%
{\it Proof.}
The sum in the first equation counts the number of times the elements of $S$ appear in the blocks of $D$ counting repetitions. Since each element appears in $D$ $r$ times, all the $S$ elements appear $r|S|$ times. Let us now count the number of times each pair of elements appears in the blocks of $D$ including all the repetitions. Every time a block of $D$ has exactly $j$ elements from $S$, this block contains ${j \choose 2}$ pairs of elements of $S$. Since each pair of elements appears $\lambda$ times in $D$ and there are ${|S|\choose 2}$ pairs in $S$, we have 
\begin{equation}
\label{choose}
\sum_{j=0}^{ k}  {j \choose 2}\delta_j(S, D) = {|S| \choose 2}  \lambda,
\end{equation} 
where  we assume ${m \choose n}=0$ if $n>m.$ The second equation of the lemma follows from \eqref{choose} and the first equation of the lemma.
$\Box$
\begin{theorem}
\label{thm_diag}
Let $M=M(D_1, D_2)$ be the mutual incidence matrix of two block designs $D_1$ and $D_2$ built on a finite set $V$. Then the diagonal entries of $MM^T$ are
\begin{equation*}
(MM^T)_{nn} = k_1(\lambda_2 k_1 -\lambda_2 +r_2), \quad\text{for}\qquad n=1, \dots, b_1.
\end{equation*} 
\end{theorem}
{\it Proof.} By definition, the $n$th row of $M$ is given by
\begin{equation*}
M_{nm} = |B_1^n \cap B_2^m| \qquad \text{for} \qquad m=1, \dots, b_2. 
\end{equation*} 
Thus the $n$th element on the diagonal of $MM^T$ is thus
\begin{equation}
\label{sum}
(MM^T)_{nn} = \sum_{m=1}^{b_2} |B_1^n \cap B_2^m|^2, 
\end{equation} 
which we can rewrite using the quantity $\delta_j(B_1^n, D_2)$ from Lemma \ref{lemma_formulas}. Namely, we regroup the terms in the  sum from \eqref{sum} according to the value of $|B_1^n \cap B_2^m|$: there are $\delta_j(B_1^n, D_2)$ terms equal to $j^2$ in the sum \eqref{sum} and thus
\begin{equation*}
(MM^T)_{nn} = \sum_{j=1}^{k_2} j^2 \delta_j(B_1^n, D_2)
\end{equation*} 
which is equal to $(\lambda_2k_1-\lambda_2+r_2)k_1$ according to the second equation of Lemma \ref{lemma_formulas} applied to $S=B_1^n$ and $D=D_2.$
$\Box$

The data of a mutual incidence matrix is naturally encoded in a multigraph as follows. 
\begin{definition}
\label{def_graph}
The mutual incidence graph of two block designs $D_1$ and $D_2$ is a multigraph $\Gamma_{D_1, D_2}$ such that 
\begin{itemize}
\item its set of vertices is $\{ B_1^1\dots, B_1^{b_1};\,B_2^1\dots, B_2^{b_2} \}$; it is composed of the blocks of the two designs,
\item the vertices $B_1^n$ and $B_2^m$ are connected by $|B_1^n\cap B_2^m|=(M(D_1, D_2))_{nm}$ edges,
\item the vertices corresponding to blocks of the same design are not connected by an edge. 
\end{itemize}
\end{definition}
Note that similarly to the incidence graph of a block design, the mutual incidence graph is naturally bipartite. We obtain the incidence graph of a block design  $D_2$ if for $D_1$ we take the simplest block design  built on the set $V$, whose 1-element blocks coincide with the elements of $V$. 
In the case $D_1=D_2$, the graph $\Gamma_{D_1, D_1}$ contains two vertices for each block of $D_1$. Let us merge each such pair of vertices deleting also the $k_1$ edges connecting them. The resulting graph $\hat\Gamma_{D_1, D_1}$   can be seen as a refinement of the {\it block intersection graph} $\Gamma_{D_1}$ of design $D_1$. For definition and results on block intersection graphs see \cite{graphs-book, graphs3, graphs5, hamiltonicity, watchman, graphs, graphs2} among many other publications. More precisely, one obtains $\Gamma_{D_1}$  from $\hat\Gamma_{D_1, D_1}$by keeping only one edge from each set of edges connecting any two of the remaining vertices. In a similar way, one can obtain all {\it $S$-block intersection graphs} \cite{graphs-book} for any set $S\subset\mathbb N$ in a similar way. 

\section{Eigenvectors and eigenvalues of $MM^T$}
\label{sect_eigen}

In this Section, as before, $D_i,\;i=1,2,$ are two block designs with parameters $(v,\;b_i,\;r_i,\;k_i,\;\lambda_i )$ built on the same finite set $V$.  The next lemma gives the first eigenvalue of $MM^T$. 
\begin{lemma}
\label{lemma_1}
Let $M=M(D_1, D_2)$ be the mutual incidence matrix of two block designs $D_1$ and $D_2$ defined by \eqref{Mij}.
  Then $r_1r_2k_1k_2$ is an eigenvalue of $MM^T$ corresponding to the eigenvector $(1,\dots, 1)^T\,.$
\end{lemma}
{\it Proof.} Immediately from definition \eqref{Mij} of the incidence matrix $M$ of $D_1$ and $D_2$, we obtain
\begin{eqnarray*}
&MM^T(1,\dots, 1)^T &= M \left( \sum_{j=1}^{b_1}{M_{j1}},\dots, \sum_{j=1}^{b_1}{M_{jb_2}} \right)^T
\\
&&=M\left(\sum_{j=1}^{b_1} |B_1^j\cap B_2^1|, \dots, \sum_{j=1}^{b_1} |B_1^j\cap B_2^{b_2}| \right)^T\,.
\end{eqnarray*}
From the definition of a block design, we have $\sum_{j=1}^{b_1} |B_1^j\cap S|=r_1|S|$ for any set $S\subset V$, where $|S|$ stands for the number of elements in $S$. This is because any given element appears $r_1$ times in the blocks of $D_1$ and thus every element of $S$ will contribute $r_1$ to the above sum. 
Applying this relation twice, we have
\begin{eqnarray*}
&MM^T(1,\dots, 1)^T &= r_1k_2M \left( 1,\dots, 1 \right)^T
\\
&&=r_1k_2\left(\sum_{j=1}^{b_2} |B_1^1\cap B_2^j|, \dots, \sum_{j=1}^{b_2} |B_1^{b_1}\cap B_2^j| \right)^T=r_1k_2r_2k_1\left( 1,\dots, 1 \right)^T\,.
\end{eqnarray*}
$\Box$

The following notation will be useful in finding other eigenvalues of $MM^T.$ Let $x$ and $y$ be a fixed pair of elements of $V$ and define  vectors $\mathcal Z_{D_1}$ and $\mathcal Z_{D_2}$ corresponding to the two block designs $D_1$ and $D_2,$ respectively:  
\begin{equation}
\label{Z1}
\mathcal Z_{D_i}(x,y):=(\Z_{D_i}^1(x,y), \dots, \Z_{D_i}^{b_i}(x,y))^T\,, \quad i=1,2,
\end{equation}
where
\begin{equation}
\label{Z1_comp}
\mathcal \Z_{D_i}^j(x,y):=\left\{\begin{array}{cc}0\,, & \text{ if } x,y\in B_i^j  \text{ or if } x,y\notin B_i^j \\1\,, &  \text{ if } x\in B_i^j \text{ and } y\notin B_i^j \\-1\,, & \text{ if } x\notin B_i^j \text{ and } y\in B_i^j\end{array}\right.\,.
\end{equation}
Naturally, we have $\Z_{D_i}^j(y,x)=-\Z_{D_i}^j(x,y)\,.$
For any block design $D_i$ built on the set $V,$ let us denote by $V_{D_i}$ the vector space spanned by the set of the vectors $\Z_{D_i}(x,y)$ \eqref{Z1}, \eqref{Z1_comp} as $x$ and $y$ run through the set $V:$ 
\begin{equation}
\label{VD}
V_{D_i}:={\rm Span}\{\Z_{D_i}(x,y)\,|\, x,y\in V\}.
\end{equation}

The next lemma  shows that $\Z_{D_1}(x,y)$ is an eigenvector of $MM^T$ for any pair $(x,y)$ of elements of $V.$
\begin{lemma}
\label{lemma_2} 
Let $M$ be the mutual incidence matrix of two block designs $D_1$ and $D_2$ defined by \eqref{Mij}. Then $(r_1-\lambda_1)(r_2-\lambda_2)$ is an eigenvalue of $MM^T.$  Its corresponding eigenspace contains the space $V_{D_1}$ defined by \eqref{VD}.
\end{lemma}
{\it Proof.}  Let us fix two elements $x,y\in V$ and subdivide the blocks of $D_i$ into four sets $D_i^I, D_i^{II}, D_i^{III} $ and $D_i^{IV}$ defined by: 
\begin{itemize}
\item $D_i^I$ is the set of blocks of $D_i$ containing $x$ and $y$,
\item $D_i^{II}$ is the set of blocks of $D_i$ containing $x$ and not containing $y$,
\item $D_i^{III}$ is the set of blocks of $D_i$ containing $y$ and not containing $x$,
\item $D_i^{IV}$ is the set of blocks of $D_i$ not containing $x$ nor $y$.
\end{itemize}

Note that definition \eqref{Z1_comp} may be rewritten using the new notation as follows: 
\begin{equation}
\label{Z1_comp-2}
\mathcal \Z_{D_1}^j(x,y):=\left\{\begin{array}{cc}0\,, & \text{ if }  B_i^j \in D_i^{I}\cup D_i^{IV} \\1\,, &  \text{ if } B_i^j\in  D_i^{II} \\-1\,, & \text{ if }  B_i^j \in D_i^{III}\end{array}\right.\,.
\end{equation}

Now, we show that $MM^T\Z_{D_1}(x,y)=(r_1-\lambda_1)M\Z_{D_2}(x,y).$ Consider the product $M^T\Z_{D_1}(x,y)$. Let us denote the $j$th row of $M^T$ by $S_j$, this row corresponds to the block $B_2^j$  and is of the form $(|B_1^1\cap B_2^j|, \dots, |B_1^{b_1}\cap B_2^j|).$ Let us multiply the $j$th line of $M^T$ by $\Z_{D_1}(x,y)$ and split the obtained sum according to the four groups of blocks of $D_1$ as follows:
\begin{eqnarray*}
&&S_j\Z_{D_1}(x,y) = \sum_{n=1}^{b_1} |B_1^n\cap B_2^j|\Z_{D_1}^n(x,y)
 = \sum_{n:\, B_1^n\in D_1^I}|B_1^n\cap B_2^j|\cdot 0
 \\
&&+\sum_{n:\, B_1^n\in D_1^{II}}|B_1^n\cap B_2^j|\cdot 1
+\sum_{n:\, B_1^n\in D_1^{III}}|B_1^n\cap B_2^j|\cdot (-1)
+\sum_{n:\, B_1^n\in D_1^{IV}}|B_1^n\cap B_2^j|\cdot 0\,.
\end{eqnarray*}
This can be rewritten applying Lemma \ref{lemma_main} to the sets of blocks $Q=D_1^{II}$ and $Q=D_1^{III}$ in the form:
\begin{equation}
\label{temp}
S_j\Z_{D_1}(x,y) =\sum_{u\in B_2^j} t_{D_1^{II}}(u) - \sum_{u\in B_2^j} t_{D_1^{III}}(u).
\end{equation}
Let us fix an element $u\in B_2^j$ different from $x$ and $y$. The value $t_{D_1^{II}}(u)$ according to its definition \eqref{t-def} is equal to the number of blocks in $D_1^{II}$ that contain $u$. Due to the definition of the set of blocks $D_1^{II}$, this number is also the number of blocks of $D_1$ containing the pair $\{u,x\}$ and not containing the triple $\{x,y,u\}.$ Let us denote the number of blocks in $D_1$ containing the triple $\{x,y,u\}$ by $\alpha.$ Then $t_{D_1^{II}}(u) = \lambda_1-\alpha.$ Analogously, for the same element $u$, the number $t_{D_1^{III}}(u)$ is the number of blocks in $D_1$ containing  the pair $\{u,y\}$ and not containing the triple $\{x,y,u\}$, that is again $t_{D_1^{III}}(u)=\lambda_1-\alpha.$ Thus the contributions of $u\in B_2^j\setminus\{x,y\}$ in \eqref{temp} cancel out and we have
\begin{equation}
\label{temp1}
S_j\Z_{D_1}(x,y) =\Big(t_{D_1^{II}}(x)- t_{D_1^{III}}(x)\Big)\chi_{B_2^j}(x)+\Big(t_{D_1^{II}}(y) - t_{D_1^{III}}(y)\Big)\chi_{B_2^j}(y),
\end{equation}
where $\chi_{B_2^j}$ is the characteristic function of $B_2^j$, that is 
\begin{equation*}
\chi_{B_2^j}(z) = \left\{\begin{array}{cc}1 & \text{ if } z\in B_2^j \\0 &  \text{ if } z\notin B_2^j\end{array}\right.\,.
\end{equation*}
Knowing that $t_{D_1^{III}}(x)=t_{D_1^{II}}(y)=0$ due to the definition of $D_1^{II}$ and $D_1^{III}$, the expression in \eqref{temp1} simplifies to:
\begin{equation}
\label{temp-short}
S_j\Z_{D_1}(x,y) =t_{D_1^{II}}(x)\chi_{B_2^j}(x) - t_{D_1^{III}}(y)\chi_{B_2^j}(y).
\end{equation}
From the definition of the sets of blocks $D_1^{II}$ and $D_1^{III}$, we compute that $t_{D_1^{II}}(x)=r_1-\lambda_1 =t_{D_1^{III}}(y) $ and thus
\begin{equation*}
S_j\Z_{D_1}(x,y) =(r_1-\lambda_1)\Big(\chi_{B_2^j}(x)- \chi_{B_2^j}(y)\Big)=\left\{\begin{array}{cc}0\,, & \text{ if }  B_2^j \in D_2^{I}\cup D_2^{IV} \\1\cdot (r_1-\lambda_1)\,, &  \text{ if } B_2^j\in  D_2^{II} \\-1\cdot (r_1-\lambda_1)\,, & \text{ if }  B_2^j \in D_2^{III}\end{array}\right.\,.
\end{equation*}
The right hand side of this equality is exactly $(r_1-\lambda_1)$ multiplied by $\Z_{D_2}^j(x,y)$ in the form \eqref{Z1_comp-2} and thus, recalling that $S_j$ stands for the $j$th row of $M^T,$ we have obtained
\begin{equation}
\label{temp-final}
M^T\Z_{D_1}(x,y) =(r_1-\lambda_1)\Z_{D_2}(x,y)\,.
\end{equation}
Now, we want to multiply relation  \eqref{temp-final} by the matrix $M$ on the left. Let us compute $M\Z_{D_2}(x,y)$ in the same way as we have done for $M^T\Z_{D_1}(x,y).$ Namely, denote by $L_j$ the $j$th row of the matrix $M$. We have
\begin{eqnarray*}
&L_j\Z_{D_2}(x,y) &= \sum_{n=1}^{b_2} |B_1^j\cap B_2^n|\Z_{D_2}^n(x,y)
 = \sum_{n:\, B_2^n\in D_2^{II}}|B_1^j\cap B_2^n|
-\sum_{n:\, B_2^n\in D_2^{III}}|B_1^j\cap B_2^n|
\\
&&=\sum_{u\in B_1^j} t_{D_2^{II}}(u) - \sum_{u\in B_1^j} t_{D_2^{III}}(u)
= t_{D_2^{II}}(x)\chi_{B_1^j}(x) - t_{D_2^{III}}(y)\chi_{B_1^j}(y)\,,
\end{eqnarray*}
where as before $\chi_B(z)$ is the characteristic function of the set $B$. By definition of the sets of blocks $D_2^{II}$ and $D_2^{III},$ we have that $ t_{D_2^{II}}(x)= t_{D_2^{III}}(y)=r_2-\lambda_2$ and thus
\begin{equation}
\label{temp-final2}
L_j\Z_{D_2}(x,y) =(r_2-\lambda_2)\Big(\chi_{B_1^j}(x)- \chi_{B_1^j}(y)\Big)=(r_2-\lambda_2)\Z_{D_1}^j(x,y)\,,
\end{equation}
where the last equality is obtained by examining the four cases according to the sets of blocks $D_1^{I}, D_1^{II}, D_1^{III}, D_1^{IV}$ the block $B_1^j$ may belong to. Written in the matrix form, equation \eqref{temp-final2} becomes 
\begin{equation}
\label{MZ2}
M\Z_{D_2}(x,y)=(r_2-\lambda_2)\Z_{D_1}(x,y)
\end{equation}
which, combined with \eqref{temp-final}, implies
\begin{equation*}
MM^T\Z_{D_1}(x,y) =(r_1-\lambda_1)(r_2-\lambda_2)\Z_{D_1}(x,y)
\end{equation*}
and proves the lemma.
$\Box$

\section{Characteristic polynomial of $MM^T$}
\label{sect_all}

In Section \ref{sect_eigen} we have found two eigenvalues of the matrix $MM^T\in Mat(b_1, \mathbb Z).$ Here we prove that the remaining eigenvalues are zero and find the dimensions of the eigenspaces corresponding to the two non-vanishing eigenvalues. To this end, we apply the established results to two block designs, one of which is arbitrary and another one, denoted by $\D$, is the simplest BIBD on the set $V$. 

Without loss of generality, we identify the set $V$ with the set of integers from $1$ to $v$, that is $V=\{1, 2, \dots, v\},$
and let $\D$ be the block design on the set $V$ with parameters $(v,b,r,k,\lambda)=(v,v,1,1,0)$. The blocks of $\D$ are simply elements of $V$ and we order them in the natural way: $n$th block of $\D$ is $\{n\}$. 

Denote by $\M=M(\D, D_2)$ the mutual incidence matrix of the following two block designs: $D_1:=\D$ and an arbitrary block design $D_2$ built on the set $V$. The matrix entries of $\M$ are either $0$ or $1$, according to \eqref{Mij}:
\begin{equation*}
\M_{ij}= |i \cap B_2^j| = \left\{\begin{array}{cc}1 & \text{ if } i\in B_2^j \\0 & \text{ if } i\notin B_2^j\end{array}\right.
\end{equation*}
and it is nothing but the incidence matrix of $D_2.$

From Lemmas \ref{lemma_1} and \ref{lemma_2}, we know that $\M\M^T\in Mat(v,\mathbb Z)$ has eigenvalue $\mu_1=r_2k_2$ corresponding to the eigenvector $(1,\dots, 1)^T$ and eigenvalue $\mu_2=r_2-\lambda_2$ corresponding to the eigenspace containing the space $V_\D$ \eqref{VD} generated by the vectors  $\Z_\D(x,y)\in\mathbb R^v$ defined by \eqref{Z1} and \eqref{Z1_comp} with $x,y\in V$; their $n$th components are
\begin{equation}
\label{ZD}
\Z^n_\D(x,y) = \delta_{xn} - \delta_{yn}\,,
\end{equation}
where $\delta_{st}$ is the Kronecker delta. The space $V_\D$ is a subspace of $\mathbb R^v$ of dimension $v-1\,.$ Indeed, if $e_1, \dots, e_v$ is the canonical basis in $\mathbb R^v$ with $(e_n)_j=\delta_{jn}$, then $V_\D$ is generated by the vectors $e_n-e_{n+1}$ with $n=1, \dots, v-1$ as for any pair $(x,y)$ such that $1\leq x<y\leq v$, we have $\Z_\D(x,y) = \sum_{n=x}^{y-1} (e_n-e_{n+1})\,.$
\begin{lemma}
\label{lemma_dimension}
For any block design $D$ built on the set $V$ of cardinality $v$, the dimension of the vector space $V_{D}$ defined by \eqref{VD} is $v-1.$
\end{lemma}
{\it Proof.}  Let $\M=M(\D, D_2)$ as above.
From \eqref{temp-final} in the proof of Lemma \ref{lemma_2} we know that $\M^T\Z_\D(x,y) = \Z_{D_2}(x,y)\,,$ and from \eqref{MZ2} we have that $\M\Z_{D_2}(x,y)=(r_2-\lambda_2)\Z_\D(x,y).$ This allows us to identify the matrices $\M$ and $\M^T$ with mutually inverse linear operators: 
\begin{equation}
\label{bijection}
\M^T: V_\D \to V_{D_2}, \qquad \frac{1}{r_2-\lambda_2} \M:  V_{D_2} \to V_\D
\end{equation}
and thus to conclude that ${\rm dim} V_\D={\rm dim} V_{D_2}$. Given that ${\rm dim} V_\D=v-1$ and that $D_2$ is an arbitrary block design, this proves the lemma.
$\Box$

Lemmas \ref{lemma_1} - \ref{lemma_dimension} provide most of the proof of the next theorem on the spectrum of the matrix $MM^T.$
\begin{theorem}
\label{thm}
Let $D_1$ and $D_2$ be two block designs on the set $V$ of cardinality $v$ with parameters $(v, b_i, r_i, k_i, \lambda_i), \; i=1,2$. Let $M=M(D_1, D_2)$ be their mutual incidence matrix defined by \eqref{Mij}. Then the matrix $MM^T\in Mat(b_1, \mathbb R)$ is of rank $v$ and has 
\begin{enumerate}
\item eigenvalue $\mu_1=r_1r_2k_1k_2$ with the eigenspace spanned by the vector $(1,\dots, 1)^T\,,$
\item eigenvalue $\mu_2=(r_1-\lambda_1)(r_2-\lambda_2)$ with eigenspace  $V_{D_1}$ defined by \eqref{VD} of dimension $v-1$.
\end{enumerate}
In particular, since $b_1\geq v,$ the characteristic polynomial of $MM^T$ has the form
\begin{equation*}
{\rm det} (MM^T-tI)= (-1)^{b_1}t^{b_1-v}(t-\mu_1)(t-\mu_2)^{v-1}\,.
\end{equation*}
\end{theorem}

{\it Proof.} We prove here that ${\rm rk} M\leq v$, which implies that $MM^T$ has at most $v$ non-zero eigenvalues. Given this and the results of Lemmas \ref{lemma_1} and \ref{lemma_2}, we conclude that $\mu_1$ and $\mu_2$ are the only non-zero eigenvalues. Indeed, due to Lemma \ref{lemma_2}, the eigenspace corresponding to $\mu_2$ contains the space $V_{D_1}$ \eqref{VD}. From Lemma \ref{lemma_dimension}, we know that  ${\rm dim} V_{D_1}=v-1$. If ${\rm rk} MM^T\leq v$, then the eigenspace corresponding to $\mu_2$ cannot be of dimension higher than $v-1$ and thus coincides with $V_{D_1}$. From ${\rm rk} MM^T\leq v$ and the fact that symmetric matrices are diagonalizable, we obtain ${\rm dim}\, {\rm Ker}(MM^T)=v-b_1.$

It thus remains to prove that ${\rm rk} M\leq v\,.$ To this end, let us represent all the blocks of $D_1$ and $D_2$ as vectors in $\mathbb R^v$ as follows. Consider the correspondence 
\begin{equation}
\label{Phi}
B_i^j \mapsto \Phi(B_i^j) \in \mathbb R^v,
\end{equation}
where the $n$th component of $\Phi(B_i^j)$ is 1 if $n\in B_i^j$ and $0$ otherwise. Let $\langle\cdot,\cdot\rangle$ be the Euclidean scalar product in $\mathbb R^v.$ Then we have
\begin{equation}
\label{scalar}
\langle \Phi(B_{i_1}^{j_1}), \Phi(B_{i_2}^{j_2})\rangle = |B_{i_1}^{j_1}\cap B_{i_2}^{j_2}|
\end{equation}
and thus we can interpret the mutual incidence matrix $M=M(D_1, D_2)$ \eqref{Mij} of $D_1$ and $D_2$ as a part of the Gram matrix of the system of vectors $\{\Phi(B_1^1),\dots, \Phi(B_1^{b_1}), \Phi(B_2^1), \dots, \Phi(B_2^{b_2}) \}$ in $\mathbb R^v.$ Therefore ${\rm rk} M\leq v$ and ${\rm rk} MM^T\leq v\,.$ Given that we have found $v$ linearly independent eigenvectors (the vector $(1,\dots, 1)^T$ and $V_{D_1}$ are linearly independent as eigenspaces corresponding to distinct eigenvalues), ${\rm rk} MM^T= v$ and thus ${\rm rk} M= v\,.$
$\Box$
\\
\begin{corollary}
\label{corollary}
Keeping the notation of Theorem \ref{thm}, the matrix $M^TM\in Mat(b_2, \mathbb R)$ is of rank $v$ having
\begin{enumerate}
\item eigenvalue $\mu_1=r_1r_2k_1k_2$ with the eigenspace spanned by the vector $(1,\dots, 1)^T\,,$
\item eigenvalue $\mu_2=(r_1-\lambda_1)(r_2-\lambda_2)$ with eigenspace  $V_{D_2}$ defined by \eqref{VD} of dimension $v-1$.
\end{enumerate}
In particular,  the characteristic polynomial of $M^TM$ is of the form
\begin{equation*}
{\rm det} (M^TM-tI)= (-1)^{b_2}t^{b_2-v}(t-\mu_1)(t-\mu_2)^{v-1}\,.
\end{equation*}
\end{corollary}
{\it Proof.} Given that $M^T=M(D_2, D_1)$, the proof is obtained by interchanging the roles of $D_1$ and $D_2$ in the proof of Theorem \ref{thm}. 
$\Box$
\\
The next corollary is obtained immediately from the proofs of Lemmas \ref{lemma_1}, \ref{lemma_2} and from Lemma \ref{lemma_dimension} and Theorem \ref{thm}.
\begin{corollary}
Let $D$ be any block design with parameters $(v,b,r,k,\lambda)$ built on a set $V$. The matrix $M(D,D)$ defined by \eqref{Mij} is a square symmetric matrix of rank $v$ having two non-zero eigenvalues $rk$ and $r-\lambda$. The eigenspace corresponding to the eigenvalue $rk$ is spanned by the vector $(1, \dots, 1)^T$ and the eigenspace corresponding to the value $r-\lambda$ is of dimension $v-1$ and is spanned by the vectors $\Z_D(x,y)$ defined by \eqref{Z1}, \eqref{Z1_comp} where $x,y$ ran through $V$. 
\end{corollary}

We end this section with a remark on the kernel of $MM^T. $ Let us form two matrices $A_1\in Mat(v\times b_1, \mathbb R)$ having the vectors  $ \Phi(B_1^j)$ as its columns and $A_2\in Mat(b_2\times v, \mathbb R)$ having the vectors  $ \Phi(B_2^j)$ as its rows. Then due to \eqref{scalar}, $M^T=A_2A_1$ and thus the ranks of $A_1$ and $A_2$ are not smaller than the rank of $M$. Since the rank of $M$ is $v$ then ${\rm rk}A_1={\rm rk}A_2=v$. Therefore, if $b_1>v$, the matrix $A_1$ has the kernel of dimension $b_1-v$. Since the kernel of $A_1$ is a subspace of the kernel of $MM^T$ and because the dimensions of the two kernels coincide, we have ${\rm Ker} A_1={\rm Ker} MM^T$.

\section{Examples}
\label{sect_examples}
Here we give some examples of application of Theorem \ref{thm}. For a mutual incidence matrix $M$, the structure of the matrix $MM^T$ is sometimes very simple: the matrix entries take only two values: all diagonal entries are equal among themselves (according to Theorem \ref{thm_diag}) as well as all the off-diagonal ones. This structure is 
 similar to the case when $M$ is the incidence matrix of one block design, see Example \ref{example2} and Section \ref{sect_def}. Example \ref{example3} shows the case of the matrix $MM^T$ with a different structure. Note that the pairs of block designs considered in these examples are friends with each other in the sense of \cite{friends} as can be straightforwardly deduced from the structure of their mutual incidence matrix.

\begin{example}
\label{P3}
{\rm
Let  $D_1=\P_3$ be the Fano projective plane, that is BIBD with parameters $(v_1=b_1=7, r_1=k_1=3, \lambda_1=1)$ built on the set $V=\{1,2,\dots,7\}$. 
More precisely, we have 
\begin{eqnarray*}
&&B^1_1=(1,2,4), \quad B_1^2=(2,3,5),\quad B_1^3=(3,4,6), \quad B_1^4=(4,5,7), 
\\
&&B_1^5=(1,5,6), \quad B_1^6=(2,6,7), \quad B_1^7=(1,3,7).
\end{eqnarray*}
Let $D_2$ be the BIBD with parameters $(v_2=7,b_2=28,r_2=12,k_2=3,\lambda_2=4),$ that is
\begin{eqnarray*}
&&B_2^1=(7, 1, 2), \;\; B^2_2=(7, 1, 4), \;\;B_2^3=(7, 1, 5), \;\; B_2^4=(7, 1, 6), \;\; B_2^5=(7, 2, 3), \;\; B_2^6=(7, 2, 4), 
\\
&& B_2^7=(7, 2, 5),\; \; B_2^8=(7, 3, 4), \; \; B_2^9=(7, 3, 5), \; \; B_2^{10}=(7, 3, 6), \; \; B_2^{11}=(7, 4, 6), \; \; B_2^{12}=(7, 5, 6), 
\\ &&B_2^{13}=(1, 2, 3), \; \; B_2^{14}=(1, 2, 5), \; \; B_2^{15}=(1, 2, 6), \; \; B_2^{16}=(1, 3, 4), \; \; B_2^{17}=(1, 3, 5), \; \; B_2^{18}=(1, 3, 6),
\\&& B_2^{19}=(1, 4, 5), \; \; B_2^{20}=(1, 4, 6), \; \; B_2^{21}=(2, 3, 4), \; \; B_2^{22}=(2, 3, 6), \; \; B_2^{23}=(2, 4, 5), \; \; B_2^{24}=(2, 4, 6),   
\\&& B_2^{25}=(2, 5, 6), \; \; B_2^{26}=(3, 4, 5), \; \; B_2^{27}=(3, 5, 6), \; \; B_2^{28}=(4, 5, 6).  
\end{eqnarray*}
Their mutual incidence matrix $M_a=M(D_1, D_2)$ is the $7\times 28$ matrix
\begin{equation*}
M_a=\left(\begin{array}{cccccccccccccccccccccccccccc}
2& 2& 1& 1& 1& 2& 1& 1& 0& 0& 1& 0& 2& 2& 2& 2& 1& 1& 2& 2& 2& 1&2& 2& 1& 1& 0& 1  
\\       
1& 0& 1& 0& 2& 1& 2& 1& 2& 1& 0& 1& 2& 2& 1& 1& 2& 1& 1& 0& 2& 2&2& 1& 2& 2& 2& 1
\\
0& 1& 0& 1& 1& 1& 0& 2& 1& 2& 2& 1& 1& 0& 1& 2& 1& 2& 1& 2& 2& 2&1& 2& 1& 2& 2& 2
\\
1& 2& 2& 1& 1& 2& 2& 2& 2& 1& 2& 2& 0& 1& 0& 1& 1& 0& 2& 1& 1& 0&2& 1& 1& 2& 1& 2
\\
1& 1& 2& 2& 0& 0& 1& 0& 1& 1& 1& 2& 1& 2& 2& 1& 2& 2& 2& 2& 0& 1&1& 1& 2& 1& 2& 2
\\
2& 1& 1& 2& 2& 2& 2& 1& 1& 2& 2& 2& 1& 1& 2& 0& 0& 1& 0& 1& 1& 2&1& 2& 2& 0& 1& 1
\\
2& 2& 2& 2& 2& 1& 1& 2& 2& 2& 1& 1& 2& 1& 1& 2& 2& 2& 1& 1& 1& 1&0& 0& 0& 1& 1& 0\end{array}\right)
\end{equation*}
and the $7\times 7$ matrix $M_aM_a^T$ is such that its coefficients take only two distinct values: 
\begin{equation*}
(M_aM_a^T)_{jj}=60 \quad\text{and}\qquad (M_aM_a^T)_{ij}= 44 \quad\text{if}\quad i\neq j.
\end{equation*}
It is easy to check that $M_aM_a^T$ has the eigenvector $(1,\dots, 1)^T$ corresponding to the eigenvalue $r_1r_2k_1k_2=324$ and  the eigenvectors $\Z_{D_1}(x,y)$ with $x,y\in \{1, \dots, 7\}$ spanning a 6-dimensional eigenspace corresponding to the eigenvalue $(r_1-\lambda_1)(r_2-\lambda_2)=16.$ Due to the values of the parameters $k_1=r_1=3$ and $\lambda_1=1,$ the vectors $\Z_{D_1}(x,y)$ have the form $e_{i_1}+e_{i_2}-e_{i_3}-e_{i_4}$ with distinct indices $i_1, i_2, i_3, i_4$ where $e_j$ are the vectors forming the canonical basis in $\mathbb R^7.$ For example, we have
\begin{equation*}
\Z_{D_1}(1,2) = \left(\begin{array}{ccccccc}0 & -1 & 0 & 0 & 1 & -1 & 1\end{array}\right)^T
\qquad\text{and}\qquad
\Z_{D_1}(1,3) = \left(\begin{array}{ccccccc}1 & -1 & -1 & 0 & 1 & 0 & 0\end{array}\right)^T.
\end{equation*}
Note that the same vectors $\Z_{D_1}(x,y)$ span an eigenspace of the $7\times 7$ matrix $M_bM_b^T$ with $M_b=M(D_1, D_1)=M(\P_3, \P_3)$ whose coefficients are
\begin{equation*}
(M_b)_{jj}=3 \quad\text{and}\qquad (M_b)_{ij}= 1 \quad\text{if}\quad i\neq j
\end{equation*}
and thus for $M_b=M(D_1, D_1)$
\begin{equation*}
(M_bM_b^T)_{jj}=15 \quad\text{and}\qquad (M_bM_b^T)_{ij}= 11 \quad\text{if}\quad i\neq j.
\end{equation*}
The eigenvalue corresponding to the 6-dimensional eigenspace spanned by vectors $\Z_{D_1}(x,y), \; x,y\in\{1,\dots, 7\}$ is $\mu_2=(r_1-\lambda_1)^2=4$. The other non-zero eigenvalue corresponding to the eigenspace spanned by $(1,\dots, 1)^T$ is $\mu_1=r_1^2k_1^2=81$  in accordance with Theorem \ref{thm}.
}
\end{example}

\begin{example}
\label{example2}
\rm{
Let $D_1$ be the simplest block design $\D$ from Section \ref{sect_all} built on the set $\{1, \dots, 7\}$, that is the block design with parameters $(v_1=b_1=7, \,r_1=k_1=1, \,\lambda_1=0)$ and $D_2$ be the Fano plane $\P_3$ from Example \ref{P3}, that is $(v_2=b_2=7, r_2=k_2=3, \lambda_2=1).$ The mutual incidence matrix $M_{\P_3}=M(\D, \P_3)$ is the usual incidence matrix of $\P_3.$ Thus $M_{\P_3}M_{\P_3}^T$ is the $7\times 7$ matrix with entries
\begin{equation*}
(M_{\P_3}M_{\P_3}^T)_{jj}=r_2=3 \quad\text{and}\qquad (M_{\P_3}M_{\P_3}^T)_{ij}= \lambda_2=1 \quad\text{if}\quad i\neq j.
\end{equation*}
As is easy to check, the rank of this matrix is $7$ and the vector $(1, \dots, 1)^T$ is an eigenvector corresponding to eigenvalue $\mu_1=9.$ The eigenvalue $\mu_2=2$ corresponds to the eigenspace of dimension 6 spanned by vectors $e_j-e_{j+1}$ with $j=1, \dots, 6$ or by the vectors $\Z_{D_1}(x,y)$ with $1\leq x,y\leq 7$ which in this case have the form $\Z_{D_1}(i,j) =e_i-e_j\,.$
}
\end{example}

\begin{example}
\label{example3}
\rm{
Let now $D_1$ be the block design with parameters $(v_1=6, \,b_1=10, \,r_1=5, \,k_1=3, \,\lambda_1=2)$ having the following blocks
\begin{eqnarray*}
 B_1^1=(6, 1, 4),\quad B_1^2=(6, 1, 5),\quad B_1^3=(6, 2, 3), \quad B_1^4= (6, 2, 4), \quad B_1^5=(6, 3, 5), 
 \\
B_1^6=(1, 2, 3),\quad B_1^7= (1, 2, 5),\quad B_1^8=(1, 3, 4),\quad B_1^9=(2, 4, 5), \quad B_1^{10}=(3, 4, 5)
 \end{eqnarray*}
and $D_2$ be the block design with parameters $(v_2=6,\, b_2=15,\, r_2=5,\, k_2=2,\, \lambda_2=1)$ composed by taking for blocks all $2$-element subsets of the set $V=\{1, \ldots, 6\}:$
\begin{eqnarray*}
B_2^1=(6, 1),\quad B_2^2=(6, 2),\quad B_2^3=(6, 3),\quad B_2^4=(6, 4), \quad B_2^5=(6, 5), 
\\ B_2^6=(1, 2), \quad B_2^7=(1, 3), \quad B_2^8=(1, 4), \quad B_2^9=(1, 5), \quad B_2^{10}=(2, 3), 
\\ B_2^{11}=(2, 4), \quad B_2^{12}=(2, 5), \quad B_2^{13}=(3, 4), \quad B_2^{14}=(3, 5), \quad B_2^{15}=(4, 5)\,.
 \end{eqnarray*}
Then their mutual incidence matrix $M_c=M(D_1, D_2)$ is
\begin{equation*}
M_c=\left(\begin{array}{ccccccccccccccc}
2& 1& 1& 2& 1& 1& 1& 2& 1& 0& 1& 0& 1& 0& 1
\\2& 1& 1& 1& 2& 1& 1& 1& 2& 0& 0& 1& 0& 1& 1
\\1& 2& 2& 1& 1& 1& 1& 0& 0& 2& 1& 1& 1& 1& 0
\\1& 2& 1& 2& 1& 1& 0& 1& 0& 1& 2& 1& 1& 0& 1
\\1& 1& 2& 1& 2& 0& 1& 0& 1& 1& 0& 1& 1& 2& 1
\\1& 1& 1& 0& 0& 2& 2& 1& 1& 2& 1& 1& 1& 1& 0
\\1& 1& 0& 0& 1& 2& 1& 1& 2& 1& 1& 2& 0& 1& 1
\\1& 0& 1& 1& 0& 1& 2& 2& 1& 1& 1& 0& 2& 1& 1
\\0& 1& 0& 1& 1& 1& 0& 1& 1& 1& 2& 2& 1& 1& 2
\\0& 0& 1& 1& 1& 0& 1& 1& 1& 1& 1& 1& 2& 2& 2
\end{array}\right)
\end{equation*}
and the matrix $M_cM_c^T$ is
\begin{equation*}
M_cM_c^T=\left(\begin{array}{cccccccccc}
21& 17& 13& 17& 13& 13& 13& 17& 13& 13
\\17& 21& 13& 13& 17& 13& 17& 13& 13& 13
\\13& 13& 21& 17& 17& 17& 13& 13& 13& 13
\\17& 13& 17& 21& 13& 13& 13& 13& 17& 13
\\13& 17& 17& 13& 21& 13& 13& 13& 13& 17
\\13& 13& 17& 13& 13& 21& 17& 17& 13& 13
 \\13& 17& 13& 13& 13& 17& 21& 13& 17& 13
\\17& 13& 13& 13& 13& 17& 13& 21& 13& 17
\\13& 13& 13& 17& 13& 13& 17& 13& 21& 17
\\13& 13& 13& 13& 17& 13& 13& 17& 17& 21
\end{array}\right).
\end{equation*}
As is easy to see, this matrix has eigenvector $(1,\dots, 1)^T$ corresponding to the eigenvalue $r_1r_2k_1k_2=5\cdot 5\cdot 3\cdot 2 = 150.$  The eigenvalue $(r_1-\lambda_1)(r_2-\lambda_2)=12$ according to Theorem \ref{thm} has the eigenspace of dimension 5 spanned by the vectors $\Z_{D_1}(x,y)$ with $x,y\in \{1, \dots, 6\}.$ For example, we have
\begin{equation*}
\Z_{D_1}(1,2) = \left(\begin{array}{cccccccccc}1 & 1 & -1 & -1 & 0 & 0 & 0&1&-1&0\end{array}\right)^T
\end{equation*}
and
\begin{equation*}
\Z_{D_1}(3,4) = \left(\begin{array}{cccccccccc}-1 & 0 & 1 & -1 & 1 & 1 & 0&0&-1&0\end{array}\right)^T,
\end{equation*}
 which are indeed eigenvectors corresponding to the eigenvalue 12. 
 The kernel of $M_cM_c^T$ is 4-dimensional and is spanned by the vectors
\begin{eqnarray*}
\left(\begin{array}{cccccccccc}1 & -1 & 1 & -1 & 0 & -1 & 1 & 0 & 0 & 0\end{array}\right)^T,
\qquad \left(\begin{array}{cccccccccc}-1 & 1 & 1 & 0 & -1 & -1 & 0 & 1 & 0 & 0\end{array}\right)^T,
 \\
\left(\begin{array}{cccccccccc}1 & 0 & 2 & -2 & -1 & -1 & 0 & 0 & 1 & 0\end{array}\right)^T,
\qquad \left(\begin{array}{cccccccccc}0 & 1 & 2 & -1 & -2 & -1 & 0 & 0 & 0 & 1\end{array}\right)^T.
 \end{eqnarray*}
Interchanging $D_1$ and $D_2$, we apply Theorem \ref{thm} to the matrix $M_c^T=M(D_2, D_1)$ and obtain the $15\times 15$ matrix $M_c^TM_c$ of the form
\begin{equation*}
\small
M_c^TM_c=\left(\begin{array}{ccccccccccccccc}
14& 11& 11& 11& 11& 11& 11& 11& 11& 8& 8& 8& 8& 8& 8 
\\11& 14& 11&  11& 11& 11& 8& 8& 8& 11& 11& 11& 8& 8& 8 
\\11& 11& 14& 11& 11& 8& 11& 8&   8& 11& 8& 8& 11& 11& 8 
\\11& 11& 11& 14& 11& 8& 8& 11& 8& 8& 11& 8& 11& 8&  11 
\\11& 11& 11& 11& 14& 8& 8& 8& 11& 8& 8& 11& 8& 11& 11 
\\11& 11& 8& 8& 8& 14& 11& 11& 11& 11& 11& 11& 8& 8& 8 
 \\11& 8& 11& 8& 8& 11& 14& 11&11& 11& 8& 8& 11& 11& 8 
   \\11& 8& 8& 11& 8& 11& 11& 14& 11& 8& 11& 8& 11&  8& 11 
   \\11& 8& 8& 8& 11& 11& 11& 11& 14& 8& 8& 11& 8& 11& 11 
   \\8& 11& 11& 8& 8& 11& 11& 8& 8& 14& 11& 11& 11& 11& 8 
   \\8& 11& 8& 11& 8& 11& 8& 11& 8& 11& 14& 11& 11& 8& 11 
   \\8& 11& 8& 8& 11& 11& 8& 8& 11& 11& 11& 14& 8& 11& 11 
   \\8& 8& 11& 11& 8& 8& 11& 11& 8& 11& 11& 8& 14& 11& 11 
   \\8& 8& 11&   8& 11& 8& 11& 8& 11& 11& 8& 11& 11& 14& 11 
    \\8& 8& 8& 11& 11& 8& 8& 11&11& 8& 11& 11& 11& 11& 14
\end{array}\right).
\end{equation*}
This matrix also has eigenvalues 150 and 12 with the eigenspace corresponding to the eigenvalue 12 being 5-dimensional and spanned by the vectors $\Z_{D_2}(x,y)$ with $x,y\in \{1, \dots, 6\}.$ The kernel of $M_c^TM_c$ is 9-dimensional. 
}
\end{example}

\vskip 1cm
{\bf Acknowledgments.} VS gratefully acknowledges
support from the Natural Sciences and Engineering Research Council of Canada (discovery grant) and the University of Sherbrooke.


\bigskip







\begin{thebibliography}{99}
 
 \bibitem{graphs3} Abueida, A. A., Pike, D. A., Cycle extensions of BIBD block-intersection graphs, J. Combin. Des. {\bf 21} (2013), 303-310.

 \bibitem{graphs5} Ahadi, A., Besharati, N., Mahmoodian, E. S., Mortezaeefar, M., {\it Silver block intersection graphs of Steiner 2-designs}
Graphs Combin. {\bf 29} (2013), no. 4, 735-746.

 \bibitem{graphs-book} Dewar, M., Stevens, B. {\it Ordering block designs}, CMS Books Math.
Springer, New York, 2012, xii+207 pp.

\bibitem{encyclopedia} Beth, T.,  Jungnickel, D., Lenz, H., {\it Encyclopedia of mathematics and its applications, } Cambridge University Press, 1999. Design Theory Second Edition, volume I. 

\bibitem{graphs2} Cameron, R. A., Pike, D. A. {\it Decomposable twofold triple systems with non-Hamiltonian 2-block intersection graphs}
Electron. J. Combin. {\bf 27} (2020), no. 4, Paper No. 4.57, 18 pp.
  
 \bibitem{handbook} Colbourn, C. J.,  Dinitz, J. H. (editors), {\it Handbook of combinatorial designs.} Second edition. Discrete Mathematics and its Applications (Boca Raton). Chapman \& Hall/CRC, Boca Raton, FL, 2007.
 
 \bibitem{crn} Crnkovi\'c, D.
{\it On path graphs of incidence graphs.} Math. Commun. {\bf 9} (2004), no.2, 181-188.

\bibitem{watchman} Dyer, D., Howell, J., {\it Watchman's walks of Steiner triple system block intersection graphs}
Australas. J. Combin. {\bf 68} (2017), 23-34.
 
\bibitem{igraph} Fern\'andez, B., Rukavina, S., {\it On the 2- Y -homogeneous condition of the incidence graphs of 2-designs}
Graphs Combin. {\bf 38} (2022), no. 4, Paper No. 119, 15 pp.

\bibitem{hamiltonicity} LeGrow, J. T., Pike, D. A., Poulin, J., {\it Hamiltonicity and cycle extensions in 0-block-intersection graphs of balanced incomplete block designs} Des. Codes Cryptogr. {\bf 80} (2016), no. 3, 421-433.
 
 \bibitem{nato} Haemers, W. H.
{\it Matrices for graphs, designs and codes.} Information security, coding theory and related combinatorics, 253-277.
NATO Sci. Peace Secur. Ser. D Inf. Commun. Secur., {\bf 29}, IOS Press, Amsterdam, 2011.
 
 \bibitem{graphs} Lu, X.-N. {\it Further results on existentially closed graphs arising from block designs}
Graphs Combin. 35 (2019), no. 6, 1323-1335.
 
 \bibitem{ma} Ma, J., {\it The distance signatures of the incidence graphs of affine resolvable designs,} Linear Algebra Appl. {\bf 493} (2016), 37-44.
 
 \bibitem{Rib}  Rybnikov, K. A., {\it Vvedenie v kombinatornyi analiz.} (Russian) [{\it Introduction to combinatorial analysis }] Izdat. Moskov. Univ., Moscow  
 (1972).
 
 \bibitem{R}  Ryser, H. J., {\it Combinatorial Mathematics}, The Carus Mathematical Monographs No. 14, 
Mathematical Association of America distributed by John Wiley and Sons, Inc., New York, 1963, xiv+154 pp.

 \bibitem{friends}  Shramchenko, A., Shramchenko, V.,  {\it A new relationship between block designs},  Int. J. Comb.
{\bf 2016}, Article ID 7092474 (2016).

\bibitem{subspace} Tabak, K., {\it Dual incidences and  $t$-designs in vector spaces},
J. Combin. Des. 32 (2024), no. 1, 46-52

\bibitem{domination} Tang, L., Zhou, S., Chen, J.,  {\it Domination number of incidence graphs of block designs},
Appl. Math. Comput. 363 (2019), 124600, 6 pp.

\end{thebibliography}
\end{document}